\documentclass[12pt,oneside]{amsart}
\usepackage{faktor}
\usepackage[hypertexnames=false]{hyperref}
\hypersetup{
	colorlinks=true,
	linktoc=all,
	linkcolor=black,
}

\usepackage{graphicx}
\usepackage{latexsym,amssymb,amscd}
\usepackage{amsmath}
\usepackage[all]{xy}
\usepackage{amsthm}
\usepackage{mathrsfs}
\usepackage{times,enumerate}
\usepackage[usenames,dvipsnames]{color}
\usepackage{blkarray}
\usepackage{amsmath}
\usepackage{faktor}
\usepackage{tikz}

\usetikzlibrary{arrows,decorations.pathmorphing,backgrounds,positioning,fit,matrix}
\usepackage{tikz-cd}
\usepackage{float}
\usepackage{hhline}
\usepackage[english]{babel}

\usepackage{mathtools}
\usepackage{stackengine}
\setstackEOL{\\}
\usepackage{tikz-cd}
\usetikzlibrary{arrows,decorations.pathmorphing,backgrounds,positioning,fit,petri,calc,shapes.misc,decorations.markings}
\tikzset{degil/.style={
		decoration={markings,
			mark= at position 0.5 with {
				\node[transform shape] (tempnode) {$\backslash$};
			}
		},
		postaction={decorate}
	}
}

      \usepackage{amssymb}
      \author[Afentoudls]{S. Afentoulidis-Almpanis}
      \address{Dept. of Mathematics, Bar-Ilan Univ., Ramat-Gan, 5290002, Israel}
      \email{spyridon.almpanis@biu.ac.il}
      
      \author{G. Liu}
    \address{Universit\'e de Lorraine, CNRS, IECL, F-57000 Metz, France}      
    \email{gang.liu@univ-lorraine.fr}
      
       \author{S. Mehdi}
    \address{Universit\'e de Lorraine, CNRS, IECL, F-57000 Metz, France}        
     \email{salah.mehdi@univ-lorraine.fr}
      
      \usepackage{lipsum}

      \theoremstyle{plain}
      \newtheorem{theorem}{Theorem}[section]

          \newtheorem{problem}[theorem]{Problem}
          
      \theoremstyle{definition}

      \theoremstyle{remark}

      \makeatletter
      \def\@setcopyright{}
      \def\serieslogo@{}
      \makeatother

      \let\OLDthebibliography\thebibliography
      \renewcommand\thebibliography[1]{
      	\OLDthebibliography{#1}
      	\setlength{\parskip}{8pt}
      	\setlength{\itemsep}{0.9pt plus 0.3ex}
      }
   
\makeatletter
\g@addto@macro{\endabstract}{\@setabstract}
\newcommand{\authorfootnotes}{\renewcommand\thefootnote{\@fnsymbol\c@footnote}}%
\makeatother

\begin{document}

%



 

   

   \title[Dirac cohomology and $\Theta$-correspondence for complex dual pairs]{Dirac cohomology and $\Theta$-correspondence\\ for complex dual pairs}

 \maketitle






   \begin{abstract}
     We study the behavior of Dirac cohomology under Howe's $\Theta$-correspondence in the case of complex reductive dual pairs. More precisely, if $(G_1,G_2)$ is a complex reductive dual pair with $G_1$ and $G_2$ viewed as real groups, we describe those Harish-Chandra modules $\pi_1$ of $G_1$ with nonzero Dirac cohomology whose $\Theta$-liftings $\Theta(\pi_1)$ still have nonzero Dirac cohomology. In this case, we compute explicitly the Dirac cohomology of $\Theta(\pi_1)$.
   \end{abstract}
\vspace{3mm}

   
 





\setcounter{tocdepth}{1}

\section{Introduction}

In the late 1970's, Howe introduced the $\Theta$-correspondence in a series of papers \cite{howe,howe2,howe1} to relate irreducible representations of reductive dual pairs $(G_1,G_2)$, i.e. of reductive groups $G_1$ and $G_2$ such that $G_1$ is the centralizer of $G_2$, and vice versa, when they are viewed as subgroups of a larger group. The $\Theta$-correspondence, or $\Theta$-lifting, and its applications have been broadly studied. In particular, the $\Theta$-correspondence has been established in many cases, including the case of Archimedean local fields \cite{howe1}, p-adic local fields with $p$ odd \cite{wald}, finite fields \cite{aubert}, dual pairs of general linear groups for arbitrary residue characteristic \cite{ming}, for orthogonal-symplectic or unitary dual pairs \cite{gan} or for quaternionic dual pairs \cite{gan1}. Our results involve the explicit description by Adams and Barbasch of the $\Theta$-correspondence in the case of complex Lie groups viewed as real groups \cite{barbadams}.

On the other hand, in the early 1970's, geometric Dirac operators were used by Parthasarathy to realize explicitly discrete series representations of a real semisimple Lie group $G$ in terms of $L^2$-sections of twisted vector bundles over the Riemannian symmetric space $G/K$ \cite{parthaThesis}. (Here $K$ is a maximal compact subgroup of $G$.)  Later, as a byproduct, he deduced a necessary condition for unitarity of highest weight modules in terms of an inequality on the infinitesimal characters, known as the Parthasarathy - Dirac inequality \cite{parthaineq}. Though this inequality provides an efficient tool to detect non-unitary modules, it does not single out their infinitesimal characters. To solve this issue, in the late 1990's, Vogan introduced an algebraic analogue of Parthasarthy's Dirac operator and defined the Dirac cohomology $H_D(X)$ for an irreducible $(\mathfrak{g},K)$-module $X$. (Here ${\mathfrak g}$ stands for the complexified Lie algebra of $G$.) Then he conjectured that nonzero Dirac cohomology reveals the infinitesimal character of a unitary irreducible $X$ \cite{vogantalks}. Vogan's conjecture was proved by Huang and Pand\v zi\'c in 2001 \cite{huangpandzic} and was extended by Kostant to the case of reductive homogeneous spaces $G/H$ \cite{K3} (known as Kostant cubic Dirac cohomology), paving the way for tremendous research activities. For finite-dimensional module Dirac cohomology coincides with the kernel of the Dirac operator and is well-understood \cite{Kostant-1999,kang-pandzic,Mehdi-2014}. We mention that the
kernel of non-cubic
Dirac operators for
finite-dimensional
modules was studied
by the first author in \cite{afentoulidis}. However, for infinite-dimensional modules, except for unitary modules, Dirac cohomology need not be the kernel of the Dirac operator and may be hard to compute. Nevertheless, just to mention a few cases, Kostant constructed quotients of Verma modules with non-zero (cubic) Dirac cohomology \cite{K3}, while Barbasch and Pand\v zi\'c computed Dirac cohomology for unipotent representations of complex groups \cite{BP2}. 

It turns out that Dirac cohomology provides an interesting invariant which can be used to classify $({\mathfrak g},K)$-modules. Indeed, modules having nonzero Dirac cohomology include most of the $A_{\mathfrak q}(\lambda)$-modules, in particular the discrete series representations, finite-dimensional modules, unitary highest weight modules, and many unipotent representations. Dirac cohomology can also be related to coherent families and characters. In \cite{MPV}, Mehdi, Pand\v zi\'c and Vogan defined a notion of algebraic Dirac index and showed that this behaves well with respect to translation functors defined for coherent families of $(\mathfrak{g},K)$-modules. Later, Mehdi, Pand\v zi\'c, Vogan and Zierau defined a Dirac index polynomial over the Cartan subalgebra of $\mathfrak{g}$ and connected this polynomial with nilpotent orbits and associated cycles of Harish-Chandra modules \cite{MPVZ}. Dirac cohomology has also become a useful tool in $C^*$-algebras and non-commutative geometry (see for instance \cite{clare1,clare2}).

In this paper, we study the behavior of the Dirac cohomology under the $\Theta$-correspondence. More precisely, let 
$(G_1, G_2)$ be a reductive dual pair (see Sections \ref{correspondence} and \ref{diraccoho} for definitions). Let 
$\pi_1$ be a Harish-Chandra module of $G_1$ and $\Theta(\pi_1)$ the $\Theta$-lifting of $\pi_1$. By convention, if 
$\pi_1$ is not in the $\Theta$-correspondence, we write $\Theta(\pi_1)=0$. Write $H_D(\pi_1)$ for the Dirac cohomology of $\pi_1$. The Dirac series of $G_1$ (resp. $G_2$) is the set of Harish-Chandra modules $\pi_1$ (resp. $\pi_2$) with nonzero Dirac cohomology $H_D(\pi_1)$ (resp. $H_D(\pi_2)$). 
\begin{problem}
Find those Harish-Chandra modules $\pi_1$ in the Dirac series of $G_1$ for which the $\Theta$-lifting 
$\Theta(\pi_1)$ belongs to the Dirac series of $G_2$. Moreover, determine explicitly 
$H_D(\Theta(\pi_1))$ whenever it is nonzero.
\end{problem}
In this paper, we will address this problem in the case when $G$ is a classical complex Lie group viewed as a real group. For complex Lie groups, to our knowledge, the study of Dirac series was initiated by Barbasch and Pand\v zi\'c in \cite{BP2}. The authors provide sufficient conditions for certain families of unitary representations of $G$ to have nonzero Dirac cohomology (see Section \ref{diraccoho}). They also conjectured necessary conditions for a representation of $G$ to be in the Dirac series. The first results supporting this conjecture were obtained in \cite{dong2,dong1}, while the full conjecture was proved recently by Barbasch, Dong and Wong in \cite{dongwong}. Our study involves some of these results. 

Our paper is organized as follows. In Section \ref{sectionzhelo}, we present an overview of the classification of irreducible representations of complex Lie groups formulated by Zhelobenko in \cite{zhelobenko}. The main results of this section are gathered in Theorem \ref{zhelo}. In Section \ref{correspondence}, we recall basics facts about the $\Theta$ - correspondence and collect the results of Barbasch and Adams \cite{barbadams} we will need. Next, in Section \ref{diraccoho}, we discuss some features of the Dirac cohomology, in particular several useful results of Dong and Wong \cite{dongwong}, as well as a list of the unipotent representations with nonzero Dirac cohomology from \cite{BP2}. Finally, Sections \ref{type2} and \ref{types1} contain our main results, i.e. an explicit study of the Dirac cohomology $H_D(\Theta(\pi_1))$ of the $\Theta$-lifting $\Theta(\pi_1)$ of $\pi_1$, when $\pi_1$ is in the Dirac series of $G_1$ for reductive dual pairs of type II and I (Theorem \ref{thmx}, Theorem \ref{thmy} respectively). 


\section{Irreducible $(\mathfrak{g},K)$-modules of complex Lie groups}\label{sectionzhelo}
In this section, we recall the classification of irreducible $(\mathfrak{g},K)$-modules given by Zhelobenko in \cite{zhelobenko} when $G$ is a complex semisimple Lie group viewed as a real group. Let $G$ be such a Lie group and $\mathfrak{g}_0$ be its Lie algebra. We fix a Cartan subalgebra $\mathfrak{h}_0$ in $\mathfrak{g}_0$ and a positive system $\Delta^+$ for the set $\Delta:=\Delta(\mathfrak{g}_0,\mathfrak{h}_0)$ of ${\mathfrak t}_0$-roots in ${\mathfrak g}_0$, with Weyl group $W(\mathfrak{g}_0,\mathfrak{h}_0)$. For every positive root $\alpha$, we choose a standard $\mathfrak{sl}_2$-triplet $\{H_\alpha,X_\alpha,X_{-\alpha}\}$ such that there exist real numbers $N_{\alpha,\beta}$ with
\begin{align*}
[X_\alpha,X_\beta]&=N_{\alpha,\beta}X_{\alpha+\beta}\\
\text{and }\hspace{2mm}N_{-\alpha,-\beta}&=N_{\alpha,\beta}\;\;\forall\beta\in\Delta.
\end{align*}
Denote by $\overline{\cdot}$ the conjugation of $\mathfrak{g}_0$ with respect to the real form of $\mathfrak{g}_0$ defined by $\sum\limits_{\alpha\in\Delta^+}\mathbb{R}\{H_\alpha,X_\alpha,X_{-\alpha}\}$. Define also the antiautomorphism
\begin{equation*}
	(\cdot)^t:\mathfrak{g}_0\rightarrow \mathfrak{g}_0
\end{equation*}
which fixes the elements of $\mathfrak{h}_0$ while, for every $\beta\in\Delta$, it maps the root vector $X_\beta$ to $X_{-\beta}$. Consider the following (real) Lie subalgebras:

\begin{align*}
	\mathfrak{k}_0:&=\{X\in\mathfrak{g}_0\mid X+\bar{X}^t=0\}\\
	&=\sum\limits_{\beta\in\Delta}i\mathbb{R}\{H_\beta\}\oplus \mathbb{R}\{X_\beta-X_{-\beta}\}\oplus i\mathbb{R}\{X_\beta+X_{-\beta}\}\\
	\mathfrak{p}_0:&=\{X\in\mathfrak{g}_0\mid X=\bar{X}^t\}\\
	&=\sum\limits_{\beta\in\Delta}\mathbb{R}\{H_\beta\}\oplus i\mathbb{R}\{X_\alpha-X_{-\beta}\}\oplus \mathbb{R}\{X_\beta+X_{-\beta}\}\\
	\mathfrak{t}_0:&=\mathfrak{k}_0\cap\mathfrak{h}_0=\sum\limits_{\beta\in\Delta} i\mathbb{R}\{H_\beta\}\\
	\mathfrak{a}_0:&=\mathfrak{p}_0\cap \mathfrak{h}_0=\sum\limits_{\beta\in\Delta} \mathbb{R}\{H_\beta\}\\
	\end{align*}
Then, as real Lie algebras, $\mathfrak{g}_0$ and $\mathfrak{h}_0$ can be decomposed into
\begin{align*}
	\mathfrak{g}_0&=\mathfrak{k}_0\oplus\mathfrak{p}_0\\
	\mathfrak{h}_0&=\mathfrak{t}_0\oplus\mathfrak{a}_0.
\end{align*}

Let us now consider the complexification $\mathfrak{g}:=(\mathfrak{g}_0)_\mathbb{C}$ of $\mathfrak{g}_0$. Let $i$ and $j$ be the imaginary units coming from the complex structure of $\mathfrak{g}_0$ and the complexification of $\mathfrak{g}_0$, respectively. As a complex Lie algebra, $\mathfrak{g}$ is isomorphic to $\mathfrak{g}_0\times \mathfrak{g}_0$, the isomorphism being given by:
	\begin{equation*}
		\mathfrak{g}=\mathfrak{g}_0\oplus j\mathfrak{g}_0\ni X_1+jX_2\mapsto (X_1+iX_2,\bar{X}_1+i\bar{X}_2)\in\mathfrak{g}_0\times \mathfrak{g}_0.
	\end{equation*}

Using the above isomorphism, we obtain
\begin{align*}
	\mathfrak{k}&:=(\mathfrak{k}_0)_\mathbb{C}\cong\{(X,-X^t)\mid X\in \mathfrak{k}_0\oplus i\mathfrak{k}_0=\mathfrak{g}_0\}\\
	\mathfrak{h}&:=(\mathfrak{h}_0)_\mathbb{C}\cong\{(H,\bar{H})\mid H\in\mathfrak{h}_0\}\oplus i\{(H,\bar{H})\mid H\in\mathfrak{h}_0\}\\
	\mathfrak{t}&:=(\mathfrak{t}_0)_\mathbb{C}\cong\{(H,-H)\mid H\in\mathfrak{t}_0\oplus i\mathfrak{t}_0=\mathfrak{h}\}\\	
	\mathfrak{a}&:=(\mathfrak{a})_\mathbb{C}\cong\{(H,H)\mid H\in\mathfrak{a}_0\oplus i\mathfrak{a}_0=\mathfrak{h}\}.
\end{align*}
Therefore, the vector dual $\mathfrak{h}^*$ of $\mathfrak{h}$ can be described in two different ways. Either as the direct sum $\mathfrak{h}^*=\mathfrak{k}^*\oplus\mathfrak{a}^*$ or as the direct product $\mathfrak{h}\cong \mathfrak{h}_0^*\times\mathfrak{h}_0^*$. In other words, an element $\lambda\in\mathfrak{h}^*$ can be described as a pair $(\mu,\nu)$ where $\mu\in\mathfrak{t}^*$, $\nu\in\mathfrak{a}^*$ and $\lambda=\mu\oplus\nu$, or as a pair $(\lambda_1,\lambda_2)$ where $\lambda_1,\lambda_2\in\mathfrak{h}_0^*$ and 
$\lambda(H_1+jH_2)=\lambda_1(H_1+iH_2)+\lambda_2(\bar{H}_1+i\bar{H}_2)$ for every $H_1,H_2\in\mathfrak{h}_0$. In particular, one has:
\begin{align*}
	\mu&=\lambda_1-\lambda_2&\lambda_1&=\frac{\mu+\nu}{2}\\
	\nu&=\lambda_1+\lambda_2&\lambda_2&=\frac{-\mu+\nu}{2}
\end{align*}

Now fix $(\lambda_1,\lambda_2)\in\mathfrak{h}_0^*\times\mathfrak{h}_0^*$, with the corresponding $\mu$ and $\nu$ as above, and assume that $\mu=\lambda_1-\lambda_2$ is integral. Consider the character $\mathbb{C}_\mu\otimes \mathbb{C}_\nu$ of the Cartan subgroup $H:=TA$ of $G$, where $T$ and $A$ are the analytic subgroups of $G$ with Lie algebras $\mathfrak{t}_0$ and $\mathfrak{a}_0$ respectively. Then form the induced representation
\begin{equation*}
	\mathrm{Ind}_B^G(\mathbb{C}_\mu\otimes\mathbb{C}_\nu\otimes 1)
\end{equation*}
where $B$ is the Borel subgroup of $G$ containing $H$, associated with the positive system $\Delta^+$, and consider its space of $K$-finite vectors
\begin{equation*}
X(\lambda_1,\lambda_2):=\mathrm{Ind}_B^G(\mathbb{C}_\mu\otimes\mathbb{C}_\nu\otimes 1)_{K\mathrm{-finite}}.
\end{equation*}
 
 \begin{theorem}[\cite{zhelobenko}] The $(\mathfrak{g},K)$-module $X(\lambda_1,\lambda_2)$ contains the $K$-type of extremal weight $\mu=\lambda_1-\lambda_2$ with multiplicity $1$. 
 	\end{theorem}
 
 	Let $L(\lambda_1,\lambda_2)$ be the unique irreducible subquotient of $X(\lambda_1,\lambda_2)$ containing the $K$-type of extremal weight $\mu$. The classification of the irreducible $(\mathfrak{g},K)$-modules is given by the following theorem.

 	\begin{theorem}[\cite{zhelobenko}]\label{zhelo}  Let $X(\lambda_1,\lambda_2)$ and $L(\lambda_1,\lambda_2)$ be as above.
 		\begin{itemize}
 		\item[1.] Every irreducible admissible $(\mathfrak{g},K)$-module is of the form $L(\lambda_1,\lambda_2)$ for some $\lambda_1,\lambda_2\in\mathfrak{h}_0^*\times\mathfrak{h}_0^*$.
 		\item[2.] Two modules $L(\lambda_1,\lambda_2)$ and $L(\lambda_1',\lambda_2')$ are equivalent if and only if there is $w\in W(\mathfrak{g}_0,\mathfrak{h}_0)$ such that $\lambda_1=w\lambda_1'$ and $\lambda_2=w\lambda_2'$.
 		\item[3.] The module $L(\lambda_1,\lambda_2)$ admits a nondegenerate hermitian form if and only if there is $w\in W(\mathfrak{g}_0,\mathfrak{h}_0)$ such that $w\mu=\mu$ and $w\nu=-\bar{\nu}$.
 	\end{itemize}

 \end{theorem}

\section{$\Theta$-lifting for classical complex Lie groups}\label{correspondence}

In this section, we discuss the $\Theta$-correspondence, still in the case when $G$ is complex Lie groups viewed as a real group. Let $(W,\langle\cdot,\cdot\rangle)$ be a finite-dimensional complex symplectic vector space and consider the symplectic group $Sp(W)$ of $W$. Let $G_1,G_2$ be two closed subgroups of $Sp(W)$. We say that the pair $(G_1,G_2)$ forms a reductive dual pair if
\begin{itemize}
	\item[(i)] the centralizer of $G_1$ (resp. $G_2$) in $G$ is $G_2$ (resp. $G_1$), i.e. $Z_G(G_1)=G_2$ (resp $Z_G(G_2)=G_1$),
	\item[(ii)] $G_1$ and $G_2$ act reductively on $W$, i.e. $W$ decomposes into a direct sum of irreducible $G_1G_2$-representations.
\end{itemize}

If $W$ does not contain proper nontrivial $G_1G_2$-invariant subspaces, we say that the pair $(G_1,G_2)$ is irreducible. It turns out that for the complex case there are two types of irreducible reductive dual pairs, pairs of type I and pairs of type II \cite{howe}:\\
\textbf{Type I:} $G_1:=O_m(\mathbb{C})$ and $G_2:=Sp_{2n}(\mathbb{C})$,\\
		\textbf{Type II:} $G_1:=GL_m(\mathbb{C})$ and $G_2:=GL_n(\mathbb{C})$.

Consider, now, the Heisenberg group $H:=H(W)$ of the symplectic space $W$, i.e. the group $H(W):=W\times \mathbb{C}$ equipped with the group law:
\begin{equation*}
	(w_1,t_1)\cdot(w_2,t_2):=(w_1+w_2,t_1+t_2+\frac{1}{2}\langle w_1,w_2\rangle).
\end{equation*}

The group $Sp(W)$ acts on $H$ by the action:
\begin{equation*}
	g\cdot(w,t):=(gw,t),\text{ for } g\in Sp(W), (w,t)\in H
\end{equation*}
and fixes the elements of the center
\begin{equation*}
	Z:=\{(0,t)\vert t\in\mathbb{C}\}\cong \mathbb{C}
\end{equation*}
of $H$. As the following theorem indicates, the center $Z$ is of great importance for the representation theory of $H$.

\begin{theorem}[Stone-von Neumann Theorem \cite{mackey,vonneu,stone}]\label{stone} Let $\psi$ be a central character of $H$. Up to isomorphism, there is a unique irreducible unitary representation of $H$ with central character $\psi$
\end{theorem}

Fix a nontrivial central character $\psi$ of $H$ and $(\rho,\mathcal{H})$ a unitary irreducible representation of $H$ with central character $\psi$. Any other central character is of the form $\psi_t(x):=\psi(tx)$, $x\in\mathbb{C}$, for some $t\in\mathbb{C}$. Let $\rho_t$ be the corresponding irreducible representation given by Theorem \ref{stone}. Let $g\in Sp(W)$ and consider the representation $\rho_t(g\cdot)$ of $H$. Since the group $Sp(W)$ acts on $H$ and leaves the center $Z$ invariant, $\rho_t(g\cdot)$ is equivalent with $\rho_t(\cdot)$. In other words, there exists an intertwining operator $\omega(g)$ between $\rho_t(g\cdot)$ and $\rho_t(\cdot)$. The map $\omega$ defines a representation of $Sp(W)$ on $\mathcal{H}$ known as the metaplectic representation of $Sp(W)$. Let $\mathrm{Irr}(\cdot)$ denotes for the irreducible representations of the corresponding group and set
\begin{align*}
	\Omega_1&:=\{\pi_1\in \mathrm{Irr}(G_1)\vert \mathrm{Hom}_{G_1}(\omega,\pi_1)\neq0\}\\
	\Omega_2&:=\{\pi_2\in \mathrm{Irr}(G_2)\vert \mathrm{Hom}_{G_2}(\omega,\pi_2)\neq0\}.
\end{align*}
Then, the conditon 
\begin{equation*}
	\{(\pi_1,\pi_2)\in \Omega_1\times \Omega_2\vert \mathrm{Hom}_{G_1G_2}(\omega,\pi_1\widehat{\otimes} \pi_2)\neq0\}
\end{equation*}
defines a bijective correspondence 
\begin{equation*}
	\pi_1\in\Omega_1\longleftrightarrow\pi_2\in\Omega_2
\end{equation*}
known as the $\Theta$-correspondence. 
Moreover, one has
\begin{equation*}
	\dim\mathrm{Hom}_{G_1G_2}(\omega,\pi_1\widehat{\otimes} \pi_2)\leq 1
\end{equation*}
for any $(\pi_1,\pi_2)\in\mathrm{Irr}(G_1)\times\mathrm{Irr}(G_2)$ \cite{sun}. We say that $\pi_1$ and $\pi_2$ are in 
$\Theta$-correspondence whenever one is the $\Theta$-lifting of the other, and
\begin{equation*}
	\dim\mathrm{Hom}_{G_1G_2}(\omega,\pi_1\widehat{\otimes} \pi_2)\neq0.
\end{equation*}

In \cite{barbadams}, the authors describe explicitly the $\Theta$-correspondence in the case of complex Lie groups viewed as real groups. Before we recall their main results below, we state the following theorem that we will use later on.

\begin{theorem}[\cite{barbadams,Przebinda1996}]\label{chara}
	Suppose $\pi_1$ is the $\Theta$-lifting of $\pi_2$ for a dual pair $(G_1,G_2)$. Let $(\lambda_1,\lambda_1^\prime)\in{\mathfrak h}_0^*\times{\mathfrak h}_0^*$ and $(\lambda_2,\lambda_2^\prime)$ be the infinitesimal characters of $\pi_1$ and $\pi_2$ respectively. Then for $\lambda_2$ and $\lambda_2^\prime$ we may take 
	\[\lambda_2=\lambda_1\cdot\widetilde{\lambda}\;\text{ and }\;\lambda_2^\prime=\lambda_1^\prime\cdot\widetilde{\lambda}\] 
	where $\widetilde{\lambda}$ takes values as follows:
	
	\begin{itemize}
		\item[1.] $(O_m(\mathbb{C}),Sp_{2n}(\mathbb{C}))$, $[m/2]\leq n:$
		\begin{equation*} \widetilde{\lambda}=(n-m/2,n-m/2-1,\ldots,1-\tau/2),
		\end{equation*}
		\item[2.] $(Sp_{2m}(\mathbb{C}),O_{n},\mathbb{C}))$, $m\leq[n/2]:$
		\begin{equation*}
			\widetilde{\lambda}=(n/2-m-1,n/2-m-3,\ldots,\tau/2),
		\end{equation*}
		\item[3.]$(GL_m(\mathbb{C}),GL_n(\mathbb{C}))$, $m\leq n:$
		\begin{equation*}
			\widetilde{\lambda}=\dfrac{1}{2}(n-m-1,n-m-2,\ldots,-n+m+1).
		\end{equation*}
	\end{itemize}
Here the symbol $\cdot$ stands for concatenation.
\end{theorem}

We can now state the results of Adams and Barbasch \cite{barbadams}, giving the explicit $\Theta$-correspondence for complex reductive Lie groups.

\begin{theorem}\cite[Theorem 2.8]{barbadams}\label{type1}
	For $\tau=0,1$, consider the reductive dual pair $(G_1,G_2)=(O_{2m+\tau}(\mathbb{C}),Sp_{2n}(\mathbb{C}))$. Let $\pi_1=L(\mu_1,\nu_1)$ be an irreducible representation of $G_1$. Define the integer $k=k[\mu_1]$ by writing $\mu_1=(a_1,\ldots,a_k,0,\ldots,0;\varepsilon)$ with $a_1\geq a_2\geq \ldots\geq a_k>0$. The parameter $\varepsilon$ takes values $\pm1$ according to \cite[6]{barbadams}. Write $\nu_1=(b_1,\ldots,b_m)$ and define the integer $0\leq q=q[\mu_1,\nu_1]\leq m-k$ to be the largest integer such that $2q-2+\tau, 2q-4+\tau,\ldots,\tau$ all occurring (in any order) in $\{\pm b_{k+1},\pm b_{k+2},\ldots,\pm b_m\}$. After conjugating by the stabilizer of $\mu_1$ in $W$, we may write 
	
	\begin{align*}
		\mu_1&=(\overbrace{a_1,\ldots,a_k}^{k},\overbrace{0,\ldots,0}^{m-q-k},\overbrace{0,\ldots,0}^{q},\varepsilon)\\
		\nu_1&=(\overbrace{b_1,\ldots,b_k}^{k},\overbrace{b_{k+1},\ldots,b_{m-k}}^{m-q-k},\overbrace{2q-2+\tau,2q-4+\tau,\ldots,\tau}^{q}).
	\end{align*}
	Let $\mu_1'=(a_1,\ldots,a_k)$, $\nu_1'=(b_1,\ldots,b_k)$, and $\nu_1''=(b_{k+1},\ldots,b_{m-q})$. Then, $\pi_1$ belongs to the $\Theta$-correspondence if and only if $n\geq m-eq+\frac{1-\varepsilon}{2}\tau$. In this case, the $\Theta$-lifting of $\pi_1$ is $L(\mu_2,\nu_2)$ where 
	\begin{align*}
		\mu_2=(&\mu_1',\overbrace{1\ldots,1}^{\frac{(1-\varepsilon)(2q+\tau)}{2}},0,\ldots,0)\\
		\nu_2=(&\nu_1',\overbrace{2q-1+\tau,2q-3+\tau,\ldots,2\varepsilon q+1+\varepsilon \tau}^{\frac{(1-\varepsilon)(2q+\tau)}{2}},\nu_1'',\\
		&2n-2m-\tau,2n-2m-2-\tau\ldots,-\varepsilon(2q+\tau)+2)
	\end{align*}
\end{theorem}

\begin{theorem}\cite[Theorem 2.9]{barbadams}\label{adams} Let $(G_1,G_2):=(GL_m(\mathbb{C}),GL_n(\mathbb{C}))$ be a reductive dual pair of type II with $m\leq n$. Let $\pi=L(\mu_1,\nu_1)$ be an irreducible representation of $G_1$. Write 
	\begin{equation*}
		\mu_1=(a_1,\ldots,a_k,\overbrace{0,\ldots,0}^{m-k-l},b_1,\ldots,b_l),
	\end{equation*}
	with $a_1\geq \ldots\geq a_k>0>b_1\geq\ldots\geq b_l$
	and $\nu_1=(c_1,\ldots,c_m)$. Then $\pi_1$ occurs in the $\Theta$-correspondence and $\pi_2=L(\mu_2,\nu_2)^*$ (the contragredient representation) with
	\begin{align*}
		\mu_2&=(a_1,\ldots,a_k,\overbrace{0,\ldots,0}^{m-l-k},\overbrace{0,\ldots,0}^{n-m},b_1,\ldots,b_l)\\
		\nu_2&=(c_1,\ldots,c_k,\overbrace{c_{k+1},\ldots,c_{m-l}}^{m-l-k},\overbrace{n-m-1,\ldots,-n+m+1}^{n-m},c_{m-l+1},\ldots,c_m)
	\end{align*}
\end{theorem}

\section{Dirac cohomology for complex Lie groups}\label{diraccoho}

Let $G$ be a real noncompact semisimple connected Lie group with complexified Lie algebra $\mathfrak{g}$, $\mathfrak{t}$ a Cartan subalgebra of $\mathfrak{g}$, $\langle\cdot,\cdot\rangle$ the Killing form of $\mathfrak{g}$, $K$ a maximal compact subgroup of $G$ of complexified Lie algebra $\mathfrak{k}$ associated with a Cartan involution $\theta$ of $G$, and
\begin{equation*}
	\mathfrak{g}=\mathfrak{k}\oplus\mathfrak{p}
\end{equation*}
the Cartan decomposition with respect to $\theta$. Let $\mathbf{C}(\mathfrak{p})$ be the Clifford algebra of $\mathfrak{p}$, i.e. the quotient of the tensor algebra $T(\mathfrak{p})$ of $\mathfrak{p}$ by its ideal generated by the elements $X\otimes Y+Y\otimes X-\langle X,Y\rangle$, $X,Y\in\mathfrak{p}$, and $S$ a spinor space of $\mathbf{C}(\mathfrak{p})$, i.e. a simple $\mathbf{C}(\mathfrak{p})$-module. The space $S$, equipped with the $\mathfrak{k}$-action defined via a composition map
\begin{equation}\label{map}
	\mathfrak{k}\overset{\mathrm{ad}}{\longrightarrow} \mathfrak{so}(\mathfrak{p})\overset{\varphi}{\longrightarrow} \mathbf{C}(\mathfrak{p})\longrightarrow \mathrm{End}(S),
\end{equation}
becomes a $\mathfrak{k}$-module, known as the spin module of $\mathfrak{k}$.

In the late 90's, Vogan introduced an algebraic Dirac operator as follows \cite{vogantalks}. For every $\mathfrak{g}$-module $(\pi,X)$, he defined the algebraic Dirac operator associated with $X$ to be the linear operator
\begin{equation}
	D_{{\mathfrak g},{\mathfrak k}}(X):X\otimes S\rightarrow X\otimes S
\end{equation}
given by
\begin{equation*}D_{{\mathfrak g},{\mathfrak k}}(X)=\sum_i \pi(Z_i)\otimes\gamma(Z_i)
\end{equation*}
where $\{Z_i\}$ is an orthonormal basis for ${\mathfrak p}$. Then, $D_{\mathfrak{g},\mathfrak{k}}(X)$ is $\mathfrak{k}$-equivariant and $\ker D_{\mathfrak{g},\mathfrak{k}}(X)$ is a $\mathfrak{k}$-representation. Vogan defined the Dirac cohomology of $ X $ to be the  $\widetilde {K} $-module
\begin{equation}\label{diraccohomfr}
	H_D^{{\mathfrak g},{\mathfrak k}}(X)=\frac{\ker D_{{\mathfrak g},{\mathfrak k}}(X)}{\text{Im}D_{{\mathfrak g},{\mathfrak k}}(X)\cap\ker D_{{\mathfrak g},{\mathfrak k}}(X)},
\end{equation}
where $\widetilde{K}$ denotes the spin double cover of $K$. It turns out that the algebraic Dirac operator $D_{\mathfrak{g},\mathfrak{k}}(X)$ and the Dirac cohomology $H_D^{{\mathfrak g},{\mathfrak k}}(X)$ are powerful tools to study representations of $\mathfrak{g}$. For instance, assuming that $X$ is irreducible, if $H_D^{{\mathfrak g},{\mathfrak k}}(X)$ contains a $\widetilde{K}$-module with highest weight $\beta$ then $X$ has infinitesimal character $\beta+\rho_{\mathfrak k}$. In other words, the infinitesimal character of a module can be recovered from its Dirac cohomology \cite{huangpandzic}. 

Dirac cohomology has been studied for various families of groups and modules. We are interested in the case where the group $G$ is a complex Lie group viewed as real. The set of irreducible unitary representations of $G$ with nonzero Dirac cohomology is called the Dirac series of $G$ and will be denoted by $\widehat{G}^d$.
In \cite{BP2}, Barbasch and Pand\v zi\' c provided a sufficient condition for a unitary representation 
$\pi:=L(\lambda_1,\lambda_2)$ of $G$ to be in the Dirac series.
Namely, they proved the following theorem.

\begin{theorem}[\cite{BP2}]\label{barbaschpandzic}
 Let $P = M N$ be a parabolic subgroup of $G$ and let $\Delta^+ = \Delta_\mathfrak{m}^+\cup\Delta(\mathfrak{n})$ be the corresponding system of positive roots. Let $\pi_\mathfrak{m}$ be an irreducible unitary representation of $M$ with nonzero Dirac cohomology, and let $\xi$ be a unitary character of $M$ which is dominant with respect to $\Delta^+$. Suppose that twice the infinitesimal character of $\pi=\mathrm{Ind}^G_P [(\pi_\mathfrak{m} \otimes \xi)\otimes 1]$ is regular and integral. Then $\pi$ has nonzero Dirac cohomology.
\end{theorem}
\noindent 
Moreover, they conjectured a sufficient condition for a representation of $G$ to have nonzero Dirac cohomology. The conjecture was proved by Barbasch, Dong and Wong:

\begin{theorem}[\cite{dongwong}]\label{dongwong}
	Let $G$ be a connected complex simple Lie group and $\pi\in \widehat{G}$ with regular and half-integral infinitesimal character. Then $\pi\in \widehat{G}^d$ if and only if $\pi$ is parabolically induced from a unipotent representation with nonzero Dirac cohomology, tensored with a unitary character of $M$.
\end{theorem}

Then notion of unipotent representation was introduced by Barbasch and Vogan in the 1980's \cite{barbaschvogan} and plays an important role in Representation Theory of reductive Lie groups. The general theory of unipotent representations requires tools beyond the scope of this work. However, in this paper, we only need unipotent representations with nonzero Dirac cohomology for classical complex Lie groups. It turns out that these representations are fairly simple and are explicitly described in \cite{BP2}. For the sake of completeness, we give a brief review of these representations. For more details, the reader is refered to \cite{BP2}.

\subsection{Type A}\label{typeA}
The regular half-integral unipotent representations of $GL_m(\mathbb{C})$ with nonzero Dirac cohomology  are all spherical, i.e. of the form $L(0,\nu)$ with
\begin{equation*}
	2\nu=(a-1,a-3,\ldots,b,b-1,\ldots,-b+1,-b, \ldots,-a+3,-a+1), 
\end{equation*}
where $\{a,b\}$ is a partition of $m$ with $a\equiv b+1\mod2$ (hence $m$ should be odd).

\subsection{Type B}\label{typeB}
The regular half-integral unipotent representations of $B_n$ with nonzero Dirac cohomology are all spherical with infinitesimal character of the form
\begin{equation*}
	2\lambda=2\lambda'=(-2b+1,-2b+3,\ldots,-1;-2a,-2a+2,\ldots,-2),
\end{equation*}
where $0<a\leq b$ are integers with $a+b=n$.

\subsection{Type C}\label{typeC}
There are $2$ half-integral unipotent representations of $C_n$ with nonzero Dirac cohomology; the oscillator representations $\pi_{\text{even}}$ and $\pi_{\text{odd}}$. Their infinitesimal characters are \begin{equation*}
	\lambda=\lambda'=\frac{1}{2}(2n-1,2n-3,\ldots,1)
\end{equation*}
in the case of $\pi_{\text{even}}$, and 
\begin{align*}
	\lambda=&\frac{1}{2}(2n-1,2n-3,\ldots,1)\\
	\lambda'=&\frac{1}{2}(2n-1,2n-3,\ldots,-1)
\end{align*}
in the case of $\pi_{\text{odd}}$. Their corresponding parameters are 
\begin{align*}
	\mu=&0\\
	\nu=&2\lambda
\end{align*}
and
\begin{align*}
	\mu=&(0,\ldots,0,1)\\
	\nu=&2\lambda,
\end{align*}
respectively.

\subsection{Type D}\label{typeD}
The regular half-integral representations of $D_n$ are of $2$ types. For $0< a\leq b$ integers with $n=a+b$, there are $2$ representations, one spherical $\pi_{\text{even}}$ with
\begin{equation*}
	\lambda=\lambda'=\frac{1}{2}(2a-1,2a-3,\ldots,3,1;2b-2,2b-4,\ldots,0)
\end{equation*}
and another one nonspherical $\pi_{\text{odd}}$ with
\begin{align*}
	\lambda=&\frac{1}{2}(2a-1,2a-3,\ldots,3,1;2b-2,2b-4,\ldots,0)\\
	\lambda'=&\frac{1}{2}(2a-1,2a-3,\ldots,3,-1;2b-2,2b-4,\ldots,0)
\end{align*}
Their parameters are 
\begin{align*}
	\mu&=0\\
	\nu&=2\lambda
\end{align*}
and
\begin{align*}
	\mu&=(0,\ldots,0,1)\\
	\nu&=2\lambda,
\end{align*}
respectively.

\section{Dirac cohomology and $\Theta$-lifting for dual pairs of type II}\label{type2}

In this section, the main result is the following theorem. We keep the previous notation.

\begin{theorem}\label{thmx}
Let $(G_1,G_2):=(GL_m(\mathbb{C}),GL_n(\mathbb{C}))$, $m\leq n$, be a complex reductive dual pair of type II. Let $\pi_1$ be an irreducible unitary representation in the Dirac series of $G_1$. Suppose that $\pi_1$ has parameters $\mu_1=\lambda_1-\lambda_2$ and $\nu_1=\lambda_1+\lambda_2$. Then the $\Theta$-lifting $\pi_2:=\Theta(\pi_1)$ belongs to the Dirac series of $G_2$, with
 \begin{equation*}
	H_D(\pi_2)=2^{[\frac{n}{2}]}V(2\lambda_2-\rho_2),
\end{equation*}
where $\rho_2$ is the half-sum of positive roots of $G_2$, and
\begin{align*}
\lambda_2&=\frac{1}{2}(\mu_2+\nu_2)\\
	\mu_2&=\mu_1\cdot(\overbrace{0,\ldots,0}^{n-m})\\
	\nu_2&=\nu_1\cdot(\overbrace{n-m-1,\ldots,-n+m+1}^{n-m}).
	\end{align*}
(As before the symbol $\cdot$ stands for concatenation.)
\end{theorem}

The rest of this section is devoted to the proof of this theorem.
Let $(G_1,G_2):=(GL_m(\mathbb{C}),GL_n(\mathbb{C}))$ and $\pi_1\in \widehat{G}_1^d$. Then, according to Theorem \ref{dongwong}, there is a parabolic subgroup $P_1=M_1N_1$ of $G_1$, a unipotent representation $\pi_{M_1}$ with nonzero Dirac cohomology, and a unitary character $\xi_1$ of $M_1$ such that 
\begin{equation*}
	\pi_1=\mathrm{Ind}_{P_1}^{G_1}((\pi_{M_1}\otimes \xi_1)\otimes 1).
\end{equation*}
The Levi factor $M_1$ of $P_1$ consists of a direct product 
\begin{equation*}M_1=M_{1,1}\times\ldots \times M_{1,k},
	\end{equation*}
where $M_{1,s}:=GL_{m_{1,s}}(\mathbb{C})$, $s=1,\ldots,k$.
Then $\pi_{M_1}$ is the tensor product
\begin{equation*}
	\pi_{M_1}=\pi_{M_{1,1}}\otimes \ldots\otimes \pi_{M_{1,k}},
\end{equation*}
where $\pi_{M_{1,s}}$ is a unipotent representation of $M_{1,s}$ with nonzero Dirac cohomology. According to Subsection \ref{typeA}, they are all spherical, i.e. of the form $\pi_{M_{1,s}}=L_{M_{1,s}}(0,\nu_{1,s})$ with
\begin{align*}
	2\nu_{1,s}=(&a_{1,s}-1,a_{1,s}-3,\ldots,b_{1,s},b_{1,s}-1,\ldots,-b_{1,s}+1,-b_{1,s},\\
	& \ldots,-a_{1,s}+3,-a_{1,s}+1), 
\end{align*}
where $\{a_{1,s},b_{1,s}\}$ is a partition of $m_{1,s}$ with $a_{1,s}\equiv b_{1,s}+1\mod2$ (hence $m_{1,s}$ should be odd).
On the other hand, $\xi_1$ consists of a direct product 
\begin{equation*}
	\xi_1=\xi_{1,1}\otimes \ldots\otimes \xi_{1,k}
\end{equation*}
where $\xi_{1,s}$ is a unitary character of $M_{1,s}$. The unitary characters of $M_{1,s}=GL_{m_{1,s}}(\mathbb{C})$ are central and of the form 
\begin{align*}
	\xi_{1,s}:Z(GL_{m_{1,s}}(\mathbb{C}))\cong(\mathbb{C},\times)&\longrightarrow (\mathbb{S}_1,\times)\\
	z=e^{x+iy}&\mapsto e^{ik_{1,s} y},
\end{align*}
where $k_{1,s}\in\mathbb{Z}$.
In other words, the corresponding linear form $\xi_{1,s}$ is trivial on $\mathfrak{a}$ while on $\mathfrak{t}$ is given by
\begin{equation*}
	\xi_{1,s}=\frac{k_{1,s}}{m_{1,s}}(1,\ldots,1).
\end{equation*}

Therefore, $\pi_1=L(\mu_1,\nu_1)$ with $\mu_1=(\xi_{1,1}\vert\ldots\vert\xi_{1,k})$ and $\nu_1=(\nu_{1,1}\vert\ldots\vert\nu_{1,k})$. If necessary, act by an element $w_1$ of $W(\mathfrak{gl}(m,\mathbb{C}))$ and write $\mu_1$ and $\nu_1$ in the form
	\begin{equation*}\mu_1=(a_1,\ldots,a_k,\overbrace{0,\ldots,0}^{m-k-l},b_1,\ldots,b_l),
\end{equation*}
with $a_1\geq \ldots\geq a_k>0>b_1\geq\ldots\geq b_l$, 
and $\nu_1=(c_1,\ldots,c_m)$, as in Theorem \ref{adams}. Then $\pi_2$, the $\Theta$-lifting of $\pi_1$, is $\pi_2=L(\mu_2,\nu_2)^*$ where
	\begin{align*}
	\mu_2&=(a_1,\ldots,a_k,\overbrace{0,\ldots,0}^{m-l-k},\overbrace{0,\ldots,0}^{n-m},b_1,\ldots,b_l)\\
	\nu_2&=(c_1,\ldots,c_k,\overbrace{c_{k+1},\ldots,c_{m-l}}^{m-l-k},\overbrace{n-m-1,\ldots,-n+m+1}^{n-m},c_{m-l+1},\ldots,c_m).
\end{align*}
Acting again with an element $w_2$ of $\Delta(\mathfrak{gl}(n,\mathbb{C}))$, we can write 
\begin{equation*}
	\mu_2=(\xi_{1,1}\vert\ldots\vert\xi_{1,k}\vert\overbrace{0,\ldots,0}^{n-m})
\end{equation*}
and 
\begin{equation}\label{lastparam}
	\nu_2=(\nu_{1,1}\vert\ldots\vert\nu_{1,k}\vert\overbrace{n-m-1,\ldots,-n+m+1}^{n-m}).
\end{equation}
Set $\mu'$ and $\nu'$ to be the weights of $\mathfrak{gl}(n-m,\mathbb{C})$ consisting of the last $(n-m)$-coordinates of $\mu_2$ and $\nu_2$ respectively. In other words, 
\begin{align*}
	\mu'&=\overbrace{(0,\ldots,0)}^{n-m}\\
	\nu'&=(n-m-1,\ldots,-n+m+1)=2\rho',
\end{align*}
where $\rho'$ stands for half the sum of the positive roots of $\mathfrak{gl}(n-m,\mathbb{C})$ viewed as the lower diagonal block of $\mathfrak{gl}(n,\mathbb{C})$. The parameters $\mu'$ and $\nu'$ correspond to the trivial representation of $GL_{n-m}(\mathbb{C})$. Indeed, the only $K$-type of the trivial representation is the trivial one, corresponding to $\mu'=0$, and the infinitesimal character of the trivial representation is exactly $\rho'$.
Let  $M_2:=M_1\times M'$ and $P_2:=M_2 \times N_2=(M_1\times M')\times N_2$. Then $L(\mu_2,\nu_2)$ is the unique irreducible subquotient of the representation 
\begin{equation*}
\mathrm{Ind}_{P_2}^{G_2}((\pi_{M_1}\otimes1)\otimes1)
\end{equation*}
containing the $K$-type with extremal weight $\mu_2$. Nevertheless, according to \cite{tadic,barbaschUnitdual}, this representation is irreducible so that \begin{equation*}L(\mu_2,\nu_2)=\mathrm{Ind}_{P_2}^{G_2}((\pi_{M_1}\otimes 1)\otimes 1).\end{equation*}
 Since $\pi_{M_1}$ has nonzero Dirac cohomology and $L(\mu_2,\nu_2)$ is parabolically induced by $\pi_{M_1}$, according to \cite{BP2}, $L(\mu_2,\nu_2)$ has nonzero Dirac cohomology. Since $L(\mu_2,\nu_2)$ is unitary, by \cite[10]{mehdi-zierauColl}, we deduce that $\pi_2=L(\mu_2,\nu_2)^*$ has nonzero Dirac cohomology.
 
 From \cite[Theorem 1.6, Conjecture 2.5]{dongwong}, the Dirac cohomology $H_D(\pi_1)$ of $\pi_1$ is 
 \begin{equation*}
 	H_D(\pi_1)=2^{[\frac{l_1}{2}]}V(2\lambda_1-\rho_1).
 \end{equation*}
 On the other hand, since $H_D(\pi_2)\neq 0$, from \cite[Theorem 1.6, Conjecture 2.5]{dongwong} again, we deduce that 
 \begin{equation*}
 	H_D(\pi_2)=2^{[\frac{l_2}{2}]}V(2\lambda_2-\rho_2).
 \end{equation*}

\section{Dirac cohomology and $\Theta$-lifting for dual pairs of type I}\label{types1}

The main result of this section is the following theorem.

\begin{theorem}\label{thmy}
	Let $(G_1,G_2)$ be a complex reductive dual pair of type I. Let $\pi_1$ be an irreducible unitary representation in the Dirac series of $G_1$. Suppose $\pi_1$ has parameters $\mu_1$ and $\nu_1$. Then the $\Theta$-lifting $\pi_2:=\Theta(\pi_1)$ of $\pi_1$ has a trivial Dirac cohomology, except when $(G_1,G_2)=(O_{2m}(\mathbb{C}),Sp_{2n}(\mathbb{C}))$ and the parameter $\nu_1$ is of the form
	\begin{equation*}
		\nu_1=(\nu_{1,1}\mid\ldots\mid \nu_{1,k}\mid 2a-1,2a-3,\ldots,3,1;2b-2,2b-4,\ldots,0),
	\end{equation*}
where each $\nu_{1,s}$ corresponds to a unipotent representation of some $GL$-factor and $b=n-m+1$. In this case, the parameter $\nu_2$ of $\pi_2=L(\nu_2,\mu_2)$ is of the form 
\begin{align*}
		\nu_2=(&\nu_{1,1}\mid\ldots\mid \nu_{1,k}\mid 2a-1,2a-3,\ldots,3,1;\\
		&2n-2m,2n-2m-2,\ldots,-2b+2)
\end{align*}
	\end{theorem}

In the rest of this section, we prove Theorem \ref{thmy}.
Let $(G_1,G_2)$ be a complex reductive dual pair of type I. Define
\begin{equation*}
\Theta(\widehat{G}_1^d):=\{\Theta(\pi_1)\vert \pi_1\in \widehat{G}_1^d\}
\end{equation*}

Suppose that $\pi_1\in \widehat{G}_1^d$. Then, according to Theorem \ref{dongwong}, there is a parabolic subgroup $P_1=M_1N_1$ of $G_1$ and a unipotent representation $\pi_{M_1}$ of $M_1$ with nonzero Dirac cohomology such that 
\begin{equation*}
	\pi_1=\mathrm{Ind}_{P_1}^{G_1}(\pi_{M_1}\otimes 1).
\end{equation*}
The Levi factor $M_1$ of $P_1$ consists of a direct product 
\begin{equation*}M_1=M_{1,1}\times\ldots \times M_{1,k}\times G(t),
\end{equation*}
where $M_{1,s}:=GL_{m_{1,s}}(\mathbb{C})$, $s=1,\ldots,k$, and $G(t)$ is a group of the same type as $G_1$ of rank $t$.
Then $\pi_{M_1}$ is a tensor product
\begin{equation*}
	\pi_{M_1}=\pi_{M_{1,1}}\otimes \ldots\otimes \pi_{M_{1,k}}\otimes\pi_{G(t)},
\end{equation*} 
where $\pi_{M_{1,s}}$ is a unipotent representation of $M_{1,s}$ with nonzero Dirac cohomology, and $\pi_{G(t)}$ is a unipotent representation of $G(t)$ with nonzero Dirac cohomology. For the sake of simplicity, set $M^1_1:=M_{1,1}\times\ldots \times M_{1,k}$  
and $\pi^1_1:=\pi_{M_{1,1}}\otimes \ldots\otimes \pi_{M_{1,k}}$ so that $M_1=M^1_1 \times G(t)$ and $\pi_{M_1}=\pi^1_1 \otimes\pi_{G(t)}$.

In what follows, we assume that $\lambda_1,\lambda_1'$ (or $\mu_1$ and $\nu_1$) correspond to $\pi_1$. Then using Theorems \ref{chara} and \ref{adams}, we find the corresponding parameters $(\lambda_2,\lambda_2')$ or $\mu_2,\nu_2$ of $\pi_2$ and we comment on the necessary conditions for the Dirac cohomology of $\pi_2$ to be nonzero. In some cases, these conditions are not fulfiled and the Dirac cohomology turns out to be trivial.

We remark that there is an infinite number of regular half-integral unipotent representations for $B_n$ or $D_n$ with nonzero Dirac cohomology while for $C_n$ there are only $2$. As a consequence, the non-triviality of Dirac cohomology is not in general preserved for the reductive dual pairs of Type I. In the rest of this section we will give a more detailed analysis of this fact.

\subsection{$\mathbf{G_1:=B_m, G_2:=C_n, m\leq n}$}
\begin{align*}
2\lambda_1=2\lambda_1'=(&2\lambda_{1,1};\ldots;2\lambda_{1,k}\mid-2b+1,-2b+3,\ldots,-1;\\
&-2a,-2a+2,\ldots,-2), \quad 0<a\leq b.
\end{align*}
\begin{align*}
	2\lambda_2=2\lambda_2'=(&2\lambda_{1,1};\ldots;2\lambda_{1,k}\mid-2b+1,-2b+3,\ldots,-1;\\
	&-2a,-2a+2,\ldots,-2;
	2n-2m-1,2n-2m-3,\ldots,1)
\end{align*}

The even coordinates $-2a,\ldots,-2$ correspond to a part of the Levi factor $M_1$. Nevertheless, they do not correspond to a $GL$-factor. More precisely, they do not correspond to a unipotent representation while they can not interact with the coordinates $\lambda_{1,1};\ldots;\lambda_{1,s}$ to give a unipotent representation of some $GL$-factor. On the other hand, no unipotent representation of $Sp$ involves even coordinates.

Therefore, $\Theta(\widehat{G}_1^d)\cap\widehat{G}_2^d=\emptyset$.

\subsection{$\mathbf{G_1:=C_m, G_2:=B_n, m\leq n}$, $\mathbf{\pi}_{\text{even}}$}
\begin{equation*}
	2\lambda_1=2\lambda_1'=(\lambda_{1,1};\ldots;\lambda_{1,k}\mid 2m-1,2m-3,\ldots,1).
\end{equation*}
\begin{align*}
	2\lambda_2=2\lambda_2'=(&\lambda_{1,1};\ldots;\lambda_{1,k}\\&\mid2m-1,2m-3,\ldots,1;2n-2m-1,2n-2m-3,\ldots,1)
\end{align*}

The coordinate $1$ appears twice in the coordinates $2m-1,2m-3,\ldots,1;2n -2m-1,2n-2m-3,\ldots,1$ and is not involved in a $GL$-factor. More precisely, they can not interact with coordinates neither from the part $\lambda_{1,1};\ldots;\lambda_{1,s}$ nor from the part $2m-1,2m-3,\ldots,1;2n-2m-1,2n-2m-3,\ldots,1$ and 
give a unipotent representation for some $GL$-factor. On the other hand, they can not both belong to a unipotent representation of some $O_\cdot(\mathbb{C})$.

Therefore
$\Theta(\widehat{G}_1^d)\cap\widehat{G}_2^d=\emptyset$.

\subsection{$\mathbf{G_1:=C_m, G_2:=B_n, m\leq n}$, $\mathbf{\pi}_{\text{odd}}$}
$\quad$\\
By using a similar argument as in the previous case, we can deduce that $\Theta(\widehat{G}_1^d)\cap\widehat{G}_2^d=\emptyset$.

\subsection{$\mathbf{G_1:=D_m, G_2:=C_n}$, $\mathbf{\pi}_{\text{even}}$}
\begin{align*}
\mu_1&=0\\
\nu_1&=(\nu_{1,1}\mid\ldots\mid \nu_{1,k}\mid 2a-1,2a-3,\ldots,3,1;2b-2,2b-4,\ldots,0),
\end{align*}
where $a,b\in\mathbb{N}$ with $0<a\leq b$.
Then
\begin{align*}
	\mu_2=0&\\
	\nu_2=(&\nu_{1,1}\mid\ldots\mid \nu_{1,k}\mid 2a-1,2a-3,\ldots,3,1;\\
	&2n-2m,2n-2m-2,\ldots,-2b+2).
	\end{align*}

The coordinates $2n-2m,2n-2m-2,\ldots ,-2b+2$ correspond to a unipotent representation of some $GL$-factor only if $2n-2m=2b-2$. Therefore, the Dirac cohomology is not preserved for the vast majority of cases. In the particular case where $b=n-m+1$, the parameter $\nu_2$ has the form
\begin{equation*}
	\nu_2=(\nu_{1,1}\mid\ldots\mid \nu_{1,k}\mid 2a-1,2a-3,\ldots,3,1;
	2b-2,2b-4,\ldots,-2b+2).
\end{equation*}
Then, the representation $\pi_2$ is the unique irreducible subquotient of the induced representation 
\begin{equation*}
	\mathrm{Ind}_{P_2}^{G_2}(\pi_{M_2}\otimes 1).
\end{equation*}
Here, $\pi_{M_2}$ is the representation
\begin{equation*}
\pi_{M_2}=\pi_{1,1}\otimes\ldots\otimes \pi_{1,k}\otimes\text{trivial}\otimes\pi_{\text{even}}	
\end{equation*}
of the Levi factor
\begin{equation*}
	M_2:=M_{1,1}\times\ldots\times M_{1,k}\times GL_{2b-1}(\mathbb{C})\times Sp_{2a}(\mathbb{C}).
\end{equation*}
and so it has nonzero Dirac cohomology.

\subsection{$\mathbf{G_1:=D_m, G_2:=C_n}$, $\mathbf{\pi}_{\text{odd}}$}
$\quad$\\
By using a similar argument as in the previous case, we can deduce that in the vast majority of cases, the Dirac cohomology is not preserved. The only exception is when $b=n-m+1$ where we have a similar situation as above.


\begin{thebibliography}{999999999}
\bibitem[AB95]{barbadams}J. Adams \& D. Barbasch \emph {Reductive dual pair correspondence for complex groups}, { J. Funct. Anal.}, \textbf{132}, 1-42 (1995), https://doi.org/10.1006/jfan.1995.1099

\bibitem[Afe21]{afentoulidis}S. Afentoulidis-Almpanis, \emph{Noncubic Dirac operators for finite-dimensional modules}, { J. Lie Theory}, \textbf{31}, 1113-1140 (2021)	

\bibitem[AMR96]{aubert}A.-M. Aubert, J. Michel \& R. Rouquier, \emph {Howe correspondence for reductive groups over finite fields}, { Duke Math. J.}, \textbf{83}, 353-397 (1996)

\bibitem[Bar96]{barbaschUnitdual}D. Barbasch, \emph{The unitary dual for complex classical Lie groups}, { Invent. Math.}, \textbf{96}, 103-176 (1989)

\bibitem[BDW22]{dongwong}D. Barbasch, C. Dong, \& K. Wong, \emph{Dirac series for complex classical Lie groups: a multiplicity-one theorem}, { Adv. Math.}, \textbf{403} pp. Paper No. 108370, 47 (2022), https://doi.org/10.1016/j.aim.2022.108370

\bibitem[BP11]{BP2}D. Barbasch \& P. Pand\v zi\'c, \emph{Dirac cohomology and unipotent representations of complex groups}, { Noncommutative Geometry And Global Analysis}, \textbf{546} pp. 1-22 (2011)

\bibitem[BV85]{barbaschvogan}D. Barbasch \& D. Vogan, \emph{Unipotent Representations of Complex Semisimple Groups}, { Annals Of Mathematics}, \textbf{121}, 41-110 (1985), http://www.jstor.org/stable/1971193

\bibitem[CHS22]{clare2}P. Clare, N. Higson \& Y. Song, \emph{On the Connes-Kasparov isomorphism, II: The Vogan classification of essential components in the tempered dual},  preprint (2022)

\bibitem[CHST22]{clare1}P. Clare, N. Higson, Y. Song \& X. Tang, \emph{On the Connes-Kasparov isomorphism, I: The reduced C*-algebra of a real reductive group and the K-theory of the tempered dual}, preprint (2022)

\bibitem[DD20]{dong1}J. Ding \& C. Dong, \emph{Unitary representations with Dirac cohomology: a finiteness result for complex Lie groups}, { Forum Math.}, \textbf{32}, 941-964 (2020), https://doi.org/10.1515/forum-2019-0295

\bibitem[Don15]{dong2}C. Dong, \emph{On a conjecture of Barbasch and Pand\v zi\'c}, { Comm. Algebra}. \textbf{43}, 3382-3388 (2015), https://doi.org/10.1080/00927872.2014.910795

\bibitem[GS17]{gan1}W.-T. Gan \& B. Sun, \emph{The Howe duality conjecture: quaternionic case}, { Representation Theory, Number Theory, And Invariant Theory}, \textbf{323} pp. 175-192 (2017)

\bibitem[GT16]{gan}W.-T. Gan, \& S. Takeda, \emph{A proof of the Howe duality conjecture}, { J. Amer. Math. Soc.}, \textbf{29}, 473-493 (2016), https://doi.org/10.1090/jams/839

\bibitem[HKP09]{kang-pandzic}J.-S. Huang, Y. Kang \& P. Pand\v zi\'c, \emph{Dirac cohomology of some Harish-Chandra modules}, { Transform. Groups}, \textbf{14}, 163-173 (2009)

\bibitem[How79]{howe}R. Howe, \emph{$\Theta$-series and invariant theory}, { Automorphic Forms, Representations And L-functions (Proc. Sympos. Pure Math., Oregon State Univ., Corvallis, Ore., 1977), Part 1}, pp. 275-285 (1979)

\bibitem[How89a]{howe2}R. Howe, \emph{Remarks on classical invariant theory}, { Trans. Amer. Math. Soc.}, \textbf{313}, 539-570 (1989), https://doi.org/10.2307/2001418

\bibitem[How89b]{howe1}R. Howe, \emph{Transcending classical invariant theory}, { J. Amer. Math. Soc.}, \textbf{2}, 535-552 (1989), https://doi.org/10.2307/1990942

\bibitem[HP02]{huangpandzic}J.-S. Huang, \& P. Pand\v zi\'c, \emph{Dirac cohomology, unitary representations and a proof of a conjecture of Vogan}, { J. Amer. Math. Soc.}, \textbf{15} pp. 185-202 (2002)

\bibitem[Kos99]{Kostant-1999}B. Kostant, \emph{A cubic Dirac operator and the emergence of Euler number multiplets of representations for equal rank subgroups}, { Duke Math. J.}. \textbf{100}, 447-501 (1999)

\bibitem[Kos03]{K3}B. Kostant, \emph{Dirac cohomology for the cubic Dirac operator}, { Studies In Memory Of Issai Schur (Chevaleret/Rehovot, 2000)}, \textbf{210} pp. 69-93 (2003)

\bibitem[Mac49]{mackey}G. Mackey, \emph{A theorem of Stone and von Neumann}, { Duke Math. J.}, \textbf{16} pp. 313-326 (1949)

\bibitem[M\' in08]{ming}A. Mínguez, \emph{Correspondance de Howe explicite: paires duales de type II}, { Ann. Sci. Éc. Norm. Supér. (4)}, \textbf{41}, 717-741 (2008), https://doi.org/10.24033/asens.2080

\bibitem[MPV17]{MPV}S. Mehdi, P. Pand\v zi\'c \& D. Vogan, \emph{Translation principle for Dirac index}, { Am. J. Math.}. \textbf{139}, 1465-1491 (2017)

\bibitem[MPVZ20]{MPVZ}S. Mehdi, P. Pand\v zi\'c, D. Vogan \& R. Zierau, \emph{Dirac index and associated cycles of Harish-Chandra modules}, { Adv. Math.}, \textbf{361} pp. 106917, 34 (2020), https://doi.org/10.1016/j.aim.2019.106917

\bibitem[MZ14a]{mehdi-zierauColl}S. Mehdi \& R. Zierau, \emph{Harmonic spinors on reductive homogeneous spaces}, { Developments And Retrospectives In Lie Theory}, \textbf{37} pp. 161-181 (2014)

\bibitem[MZ14b]{Mehdi-2014}S. Mehdi, \& R. Zierau, \emph{The Dirac cohomology of a finite dimensional representation}, { Proceedings Of The American Mathematical Society}, \textbf{142} (2014,2)

\bibitem[Par72]{parthaThesis}R. Parthasarathy, \emph{Dirac operator and the discrete series}, { Ann. Of Math. (2)}, \textbf{96} pp. 1-30 (1972)

\bibitem[Par80]{parthaineq}R. Parthasarathy, \emph{Criteria for the unitarizability of some highest weight modules}, {\em Proc. Indian Acad. Sci. Sect. A Math. Sci.}, \textbf{89}, 1-24 (1980)

\bibitem[Prz96]{Przebinda1996}T. Przebinda, \emph{The duality correspondence of infinitesimal characters}, {\emph Colloquium Mathematicae}, \textbf{70}, 93-102 (1996), http://eudml.org/doc/210401

\bibitem[Sto30]{stone}M. Stone, \emph{Linear transformations in Hilbert space. III: Operational methods and group theory}, {\em Proc. Natl. Acad. Sci. USA}, \textbf{16} pp. 172-175 (1930)

\bibitem[SZ20]{sun}B. Sun, \& C.-B. Zho, \emph{Local theta correspondence: the basic theory}, (2020)

\bibitem[Tad09]{tadic}M. Tadi\' c, \emph{$GL(n,\mathbb{C})$ and $GL(n,\mathbb{R})$}, {Automorphic Forms And L-functions II. Local Aspects}, \textbf{489} pp. 285-313 (2009), https://doi.org/10.1090/conm/489/09551

\bibitem[Vog97]{vogantalks}D. Vogan, \emph{Dirac operator and unitary representations}, 3 talks at MIT Lie groups seminar,1997)

\bibitem[vNe31]{vonneu}J. v. Neumann, \emph{Die Eindeutigkeit der Schrödingerschen Operatoren}, { Math. Ann.}, \textbf{104}, pp. 570-578 (1931)

\bibitem[Wal90]{wald}J. Waldspurger, \emph{D\' emonstration d'une conjecture de dualit\' e de Howe dans le cas p-adique, $p\neq 2$}, Festschrift In Honor Of I. I. Piatetski-Shapiro On The Occasion Of His Sixtieth Birthday, Part I (Ramat Aviv, 1989). \textbf{2} pp. 267-324 (1990)

\bibitem[Zhe74]{zhelobenko}D. Zhelobenko, \emph{Garmonicheskii analiz na poluprostykh kompleksnykh gruppakh Li}, (Izdat. “Nauka”, Moscow,1974)

\end{thebibliography}
\end{document}